\documentclass[12pt,a4paper,reqno]{amsart}
\pdfoutput=1
\usepackage[hidelinks,pdfstartview=FitH]{hyperref} 
\usepackage{amssymb,tikz}
\usepackage[all,cmtip]{xy}
\usepackage{verbatim}
\usepackage[USenglish]{babel}
\usepackage[top=3cm,right=2.5cm,left=2.5cm,bottom=3cm,marginparsep=15pt,
marginparwidth=2cm,footskip=20pt]{geometry}

\allowdisplaybreaks[4]

\newtheorem{proposition}{Proposition}[section]
\newtheorem{lemma}[proposition]{Lemma}
\newtheorem{theorem}[proposition]{Theorem}

\newtheorem{definition}[proposition]{Definition}

\newtheorem*{thm*}{Theorem}

\numberwithin{equation}{section}

\renewcommand{\emptyset}{\varnothing}
\newcommand{\id}{\mathrm{id}}

\title[Leavitt path algebras of trimmable graphs]{A graded pullback structure of\\ Leavitt path algebras of trimmable graphs} 

\author[P. M.~Hajac]{Piotr M.~Hajac}
\address[P. M.~Hajac]{Instytut Matematyczny, Polska Akademia Nauk, ul. \'Sniadeckich 8, Warszawa, 00-656 Poland}
\email{pmh@impan.pl }

\author[A.~Kaygun]{Atabey Kaygun} 
\address[A.~Kaygun]{Istanbul Technical University, Istanbul, Turkey.}
\email{kaygun@itu.edu.tr}

\author[M.~Tobolski]{Mariusz Tobolski}
\address[M.~Tobolski]{Instytut Matematyczny, Polska Akademia Nauk, ul. \'Sniadeckich 8, Warszawa, 00-656 Poland}
\email{mtobolski@impan.pl }


\begin{document}
\baselineskip15pt
\parskip=0.5\baselineskip

\begin{abstract}
  Motivated by recent results in graph C*-algebras concerning an equivariant pushout structure 
  of the Vaksman-Soibelman quantum odd spheres, we introduce  a class of graphs called \emph{trimmable}.
  Then we show that the
  Leavitt path algebra of a trimmable graph is graded-isomorphic to a pullback
  algebra of simpler Leavitt path algebras and their tensor products. 
  \end{abstract}
\maketitle

\section*{Introduction}\noindent

The goal of this paper is to introduce and apply the concept of
a \emph{trimmable graph}.
We begin by recalling the fundamental concepts
of path algebras~\cite{ass06} and \emph{Leavitt path
  algebras}~\cite{L62, aa05, A15, aasm17}.  Then we define a \emph{trimmable
  graph}, and prove our main result: There is a $\mathbb{Z}$-graded
algebra isomorphism from the Leavitt path algebra of a trimmable
graph to an appropriate pullback algebra. The graph C*-algebraic version of this result
is proven in~\cite{dht}, where it was used to analyze the
generators of K-groups of quantum complex projective spaces.

\section{Leavitt path algebras}\noindent

\begin{definition}
  A \emph{graph} $Q$ is a quadruple $(Q_0,Q_1,s,t)$ consisting of the
  set of vertices~$Q_0$, the set of edges $Q_1$, and the source and
  target maps $s,t\colon Q_1\to Q_0$ assigning to each edge its
  source and target vertex respectively.
\end{definition}
\noindent
We say that a graph $Q'=(Q'_0,Q'_1,s',t')$ is a sub-graph of a
graph $Q=(Q_0,Q_1,s,t)$ iff $Q'_0\subseteq Q_0\,$,
$Q'_1\subseteq Q_1\,$, and the source and target maps $s'$ and $t'$
are respective restrictions-corestrictions of the source and target
maps $s$ and~$t$. Furthermore, we say that two edges are composable
if the end of one of them is the beginning of the other.  Now we can
define a path in a graph as a sequence of composable edges. The
length of a path is the number of edges it consists of, infinity
included. We treat vertices as zero-length paths that begin and end
in themselves. 

\begin{definition}
  Let $k$ be a field and $Q$ a graph. The \emph{path algebra} $kQ$ is
  the $k$-algebra whose underlying vector space has as its basis the
  set of all finite-length paths~$\mathrm{Path}(Q)$. 
  The product is given by the composition of paths when the end of one
  path matches the beginning of the other path.  The product is
  defined to be zero otherwise.
\end{definition}
\noindent
One can check that the path algebra $kQ$ is unital if and only if the
set of vertices $Q_0$ is finite. Then the unit is the sum of all
vertices.  It is also straightforward to verify that $kQ$ is
$\mathbb{N}$-graded by the path length.

To define a Leavitt path algebra, we need ghost edges.  For any
graph $Q=(Q_0,Q_1,s,t)$, we create a new set
$Q_1^*:=\{x^*\;|\;x\in Q_1\}$ and call its elements  \emph{ghost
  edges}. Now, the source and the target maps for the extended graph
$\widehat{Q}:=\left(Q_0,Q_1\coprod Q_1^*, \hat{s},\hat{t}\,\right)$
are defined as follows:
\begin{equation}
\hat{s}(x):=s(x),\quad \hat{s}(x^*):=t(x),\quad \hat{t}(x):=t(x),\quad \hat{t}(x^*):=s(x).
\end{equation}
\begin{definition}
Let $k$ be a field and $Q$ a graph. The {\em Leavitt path algebra} $L_k(Q)$ of a graph $Q$ is the path algebra of the extended graph
$\widehat{Q}$ divided by the ideal generated by the relations:
\vspace*{-3mm}
\begin{enumerate}
\item For all edges $x_i,x_j\in Q_1\,$, we have $x_i^*x_j=\delta_{ij}t(x_i)$.
\item For every vertex $v\in Q_0$ whose preimage  $s^{-1}(v)$  is not empty and finite, we have
\[
\sum_{x\in s^{-1}(v)}xx^*=v.
\]
\end{enumerate}
\end{definition}
\noindent
In other words, the Leavitt path algebra $L_k(Q)$ of a graph $Q$ is the universal $k$-algebra generated by
the elements $v\in Q_0\,$,  $x\in Q_1\,$, $x^*\in Q_1^*$, subject to relations:
\vspace*{-3mm}
\begin{enumerate}
\item[(L1)] $v_iv_j=\delta_{ij}v_i$ for all $v_i,v_j\in Q_0$,
\item[(L2)] $s(x)x=xt(x)=x$ for all $x\in Q_1$,
\item[(L3)] $t(x)x^*=x^*s(x)=x^*$ for all $x^*\in Q_1^*$,
\item[(L4)] $x_i^*x_j=\delta_{ij}t(x_i)$ for all $x_i,x_j\in Q_1$, and 
\item[(L5)] $\sum_{x\in s^{-1}(v)}xx^*=v$  for all $v\in Q_0$ such that  $s^{-1}(v)$  is finite and nonempty.
\end{enumerate}

Furthermore, note that the $\mathbb{N}$-grading of the path algebra $k\widehat{Q}$
induces a $\mathbb{Z}$-grading of the Leavitt path algebra $L_k(Q)$ by
counting the length of any ghost edge as~$-1$ (see~\cite[Lemma~1.7]{aa05}).
Let us recall now the Graded Uniqueness Theorem \cite[Theorem~4.8]{t-m07} that shows the importance of
this grading.  We will need it in the next section.
\begin{theorem}[\cite{t-m07}]\label{gut}
  Let $Q$ be an arbitrary graph and $k$ be any field. If $A$ is a
  $\mathbb{Z}$-graded ring, and $f:L_k(E)\to A$ is a graded ring
  homomorphism with $f(v)\neq 0$ for every vertex $v\in Q_0$, then $f$
  is injective.
\end{theorem} 

Recall that a vertex $v\in Q_0$ is called a \emph{sink} if $s^{-1}(v)=\emptyset$. Next, let $x=x_1x_2\ldots x_n$
be a path in~$Q$. If the length of $x$ is at least~$1$, and if $s(x)=v=t(x)\in Q_0$,
we say that $x$ is a \emph{closed path based at $v$}.
If in addition $s(x_i)\neq s(x_j)$ for every $i\neq j$, then $x$ is called a \emph{cycle based~at~$v$}.

Next, if $v\in Q_0$ is either a sink or a base of a cycle of length 1 (a loop), then 
a singleton set $\{v\}\subseteq Q_0$   is a basic example of a hereditary
  subset of $Q_0$~\cite[Definition~2.0.5~(i)]{aasm17}, and
it follows from Corollary~{2.4.13}~{(i)} in \cite{aasm17} (cf.~\cite[Theorem~5.7~(2)]{t-m07}) that
\begin{equation}
L_k(Q)\slash I(v)\cong L_k(Q\setminus \{v\}).
\end{equation}
 Here the graph $Q\setminus \{v\}$ is obtained by removing from $Q$ the
  vertex $v$ and every edge that ends in~$v$.  In other  words, 
\begin{equation}
(Q\setminus \{v\})_0=Q_0\setminus\{v\} \quad\text{ and }\quad
    (Q\setminus\{v\})_1=\{x\in Q_1:t(x)\neq v\}.
\end{equation}
By Lemma~2.4.1 in \cite{aasm17} (cf.~\cite[Lemma~5.6]{t-m07}), we also know that
\begin{equation}\label{herideal}
  I(v)=\left\{\sum_{i=1}^nk_ix_iy^*_i\quad\Big|\quad n\geq 1,\; x_i\,,y_i\in {\rm Path}(Q),\; \hat{t}(x_i)=\hat{s}(y^*_i)=v\right\}.
\end{equation}

\section{Trimmable graphs}

We are now ready for the main definition of the paper. Merely to focus
attention, we assume henceforth that graphs are finite, i.e.\ that
the set of vertices and the set of edges are both finite.
\begin{definition}
  Let $Q$ be a finite graph consisting of a sub-graph $Q'$ emitting at least one
  edge to an external vertex $v_0$ whose only outgoing edge $x_0$
  is a loop.  We call such a graph $(Q',v_0)$-\emph{trimmable} iff 
all  edges from $Q'$ to $v_0$ begin in a vertex emitting an edge that ends
  inside  $Q'$. 
\end{definition}
\noindent 
In symbols, a trimmable graph is described as follows:
\begin{gather}
Q_0=Q'_0\cup\{v_0\},\quad v_0\not\in Q'_0\,,\quad Q_1=Q'_1\cup t^{-1}(v_0),
\\
s^{-1}(v_0)=\{x_0\},\quad t(x_0)=v_0\,, \quad t^{-1}(v_0)\setminus\{x_0\}\neq\emptyset,
\\ \label{trim}
\forall\; v\in s\big(t^{-1}(v_0)\setminus\{x_0\}\big)\,\colon\quad s^{-1}(v)\setminus t^{-1}(v_0)\neq\emptyset.
\end{gather}
The condition for a trimmable graph guarantees the fact
that when we remove the distinguised vertex $v_0$, the resulting
graph does not have new sinks.  One can imagine a
$(Q',v_0)$-trimmable graph like this:
\begin{center}
\includegraphics[angle=0,scale=0.5]{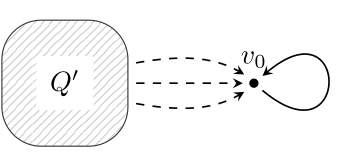}.
\end{center}

  The following graph is a simple example of a trimmable graph:
  \begin{equation}\label{trimex}
    \xymatrix{
      v_2\\
      v_1 \ar[u] \ar[r] & v_0 \ar@(u,r)
    }_{\quad .}
  \end{equation}
 Note that  the direction of the edge joining vertices $v_1$ and $v_2$ is important as the following graph is no longer trimmable:
  \begin{equation}\label{matrix}
    \xymatrix{
      v_2 \ar[d] \\
      v_1 \ar[r] & v_0 \ar@(u,r)
   }_{\quad .}
  \end{equation}

We need the trimmability condition to guarantee the existence of  maps given in the  lemma  below. This lemma is an algebraic incarnation
of a graph-C*-algebraic lemma proved in~\cite{dht}. 
The only difference between their proofs  is that instead of
using the Gauge Uniqueness Theorem \cite[Theorem~4.8]{t-m07} we use 
the Graded Uniqueness Theorem \cite[Theorem~2.3]{aHR97} (Theorem~\ref{gut}).
\begin{lemma}\label{maps}
Let $Q$ be a $(Q'_0,v_0)$-trimmable graph. Denote by $Q''$ the sub-graph of $Q$ obtained  by removing the edge~$x_0$.
The following formulas define homomorphisms of algebras:
\vspace*{-3mm}
\begin{enumerate}

\item $\pi_1\colon L_k(Q)\to L_k(Q')$,
\[ \pi_1(\alpha) = 
   \begin{cases}
     \alpha & \text{ if } \alpha\in Q'_0 \cup Q'_1\cup Q'^*_1,\\
     0 & \text{ otherwise}.
   \end{cases}
\]
\item
$\pi_2\colon L_k(Q'')\to L_k(Q')$,
\[ \pi_2(\alpha) = 
   \begin{cases}
     \alpha & \text{ if } \alpha\in Q'_0 \cup Q'_1\cup Q'^*_1,\\
     0 & \text{ otherwise}.
   \end{cases}
\]
\item
$f\colon L_k(Q)\to L_k(Q'')\otimes k[u,u^{-1}]$,
\[ f(\alpha) = 
\begin{cases}
  \alpha\otimes 1 & \text{ if } \alpha\in Q_0,\\
  v_0\otimes u & \text{ if } \alpha = x_0,\\
  v_0\otimes u^{-1} & \text{ if } \alpha = x_0^*,\\
  \alpha\otimes u & \text{ if } \alpha\in Q_1\setminus\{x_0\},\\
  \alpha\otimes u^{-1} & \text{ if } \alpha\in Q_1^*\setminus\{x_0^*\}.
\end{cases}
\]
\item
$\delta\colon L_k(Q') \to L_k(Q')\otimes k[u,u^{-1}]$,
\[ \delta(\alpha) = 
\begin{cases}
  \alpha\otimes 1 & \text{ if } \alpha\in Q_0',\\
  \alpha\otimes u & \text{ if } \alpha\in Q_1',\\
  \alpha\otimes u^{-1} & \text{ if } \alpha\in (Q_1')^*.
\end{cases}
\]
\end{enumerate}
These morphisms are $\mathbb{Z}$-graded for the standard grading on
$L_k(Q)$, $L_k(Q')$, $L_k(Q'')$, and the gradings on
$L_k(Q')\otimes k[u,u^{-1}]$ and $L_k(Q'')\otimes k[u,u^{-1}]$ 
given by the standard grading of the rightmost tensorand. Furthermore, $\pi_1$
and $\pi_2$ are surjective, and $f$ and $\delta$ are injective.
\end{lemma}

\section{A graded pullback structure}

To prove the theorem of the paper, we need a general lemma along the lines of \cite[Proposition~3.1]{pedersen}.
We omit  its routine proof.
\begin{lemma}\label{pullback}
Let $A_1$, $A_2$, $B$ and $P$ be abelian groups. A commutative diagram of group homomorphisms
{\[
\xymatrix{
&P\ar[ld]_{p_1}\ar[rd]^{p_2}&\\
A_1\ar[dr]_{q_1}&& A_2
\ar[dl]^{q_2}\\
&B&
}
\]}
is a pullback diagram if and only if the following conditions hold:
\begin{gather}
\ker(p_1)\cap\ker(p_2)=\{0\},\label{con1}\\
q_1^{-1}(q_2(A_2))=p_1(P),\label{con2}\\
p_2(\ker(p_1))=\ker(q_2).\label{con3}
\end{gather}
\end{lemma}

Recall that to prove that an algebra $P$ is a pullback algebra, one can proceed as follows. The first step is to establish
the existence of a commutative diagram of algebra homomorphisms as above. 
This implies that $p_1$ and $p_2$ define an algebra homomorphism $p$ into the pullback algebra of $A_1$ and $A_2$ over~$B$. 
Then one only needs to prove that the three conditions of Lemma~\ref{pullback} are satisfied to conclude that $p$ is an
isomorphism. Note that \eqref{con1} is equivalent to the injectivity of~$p$, whereas the conjunction of \eqref{con2} and \eqref{con3}
is equivalent to the surjectivity of~$p$.

Much as Lemma~\ref{maps}, the theorem  of the paper is an algebraic incarnation
of a graph-C*-algebraic theorem proved in~\cite{dht}. 
This time, the only difference between their proofs  is that instead of
using  \cite[Lemma~3.1]{aHR97} we use \cite[Lemma~2.4.1]{aasm17} (see \eqref{herideal}, cf.~\cite[Lemma~5.6]{t-m07}).
\begin{theorem}\label{mainthm}
Let $\pi_1$, $\pi_2$, $f$ and $\delta$ be as in Lemma~\ref{maps}.
Then the commutative diagram
{\begin{equation*}
\xymatrix{
&L_k(Q)\ar[ld]_{\pi_1}\ar[rd]^{f}&\\
L_k(Q')\ar[dr]_{\delta\quad}&& 
L_k(Q'')\otimes k[u,u^{-1}]
\ar[dl]^{\quad\pi_2\otimes\id}\\
&L_k(Q')\otimes k[u,u^{-1}]&
}
\end{equation*}}
of graded algebra homomorphisms is a pullback diagram.
\end{theorem}
\noindent
Representing pictorially the Leavitt path algebras by their respective graphs, the above diagram becomes:
\begin{center}
\vspace*{0mm}
\begin{equation*}
\xymatrix{
&\includegraphics[angle=0,scale=.3]{graphA.png}\ar[ld]_{\pi_1\quad}\ar[rd]^{\; f}&\\
\includegraphics[angle=0,scale=.3]{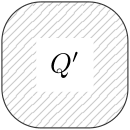}\ar[dr]_{\delta\quad}&& 
\includegraphics[angle=0,scale=.3]{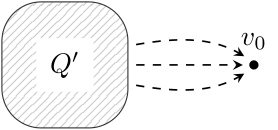}\overset{\otimes}{\phantom{\text{\huge X}}} 
\includegraphics[angle=0,scale=.26]{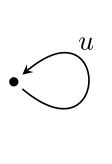}
\ar[dl]^{\quad\pi_2\otimes\id}\\
&\includegraphics[angle=0,scale=.3]{graphC.png}\overset{\otimes}{\phantom{\text{\huge X}}} 
 \includegraphics[angle=0,scale=.26]{graph0.png}.&
}
\end{equation*}
\end{center}
Note that the only non-standard map in this diagram is~$f$. It can be described verbally by  the assignment
\begin{align*}
\text{vertex}&\longmapsto\text{vertex}\otimes 1,\\
\text{$v_0$-emitted edge}&\longmapsto v_0\otimes u,\\
\text{other edge}&\longmapsto\text{other edge}\otimes u.
\end{align*}

\section*{Acknowledgement}\noindent
The authors are very grateful to Francesco D'Andrea for sharing his mathematical insight with us, and  for designing the pictures of graphs.
 It is also a pleasure to thank S\o ren Eilers, Adam Skalski and Wojciech Szyma\'nski for discussions. This
work was partially supported by NCN grant UMO-2015/19/B/ST1/03098.

\end{document}